\documentclass[a4paper]{article}

\usepackage{amsmath}
\usepackage{amssymb}
\usepackage{amsfonts}
\usepackage{enumerate,vmargin}

\newcommand{\tens}{\otimes}

\newcommand{\M}{\ensuremath{\mathbb{M}}}

\newcommand{\B}{\ensuremath{\mathbb{B}}}
\newcommand{\R}{\ensuremath{\mathbb{R}}}
\newcommand{\C}{\ensuremath{\mathbb{C}}}
\newcommand{\E}{\ensuremath{\mathbb{E}}}
\newcommand{\Id}{\ensuremath{\mathrm{Id}}}

\newcommand{\sk}{\smallskip}
\newcommand{\mk}{\medskip}

\renewcommand{\leq}{\ensuremath{\leqslant}}
\renewcommand{\geq}{\ensuremath{\geqslant}}

\newcommand{\la}{\langle}
\newcommand{\ra}{\rangle}
\newcommand{\qed}{\hfill \vrule height6pt  width6pt depth0pt}
\newtheorem{thm}{Theorem}
\newtheorem{defi}[thm]{Definition}
\newtheorem{prop}[thm]{Proposition}
\newtheorem{cor}[thm]{Corollary}

\newenvironment{pf}[1][]{\noindent {\it Proof #1} : }{\hbox{~}\qed
\smallskip
}

\title{A Markov dilation for self-adjoint Schur multipliers}

\date{}
\author{\'Eric Ricard}
\begin{document}

\maketitle

\begin{abstract}

 We give a formula for Markov dilation in the sense of Anantharaman-Delaroche
for real positive Schur multipliers on $\B(H)$.

\end{abstract}

\makeatletter
\renewcommand{\@makefntext}[1]{#1}
\makeatother
\footnotetext{\noindent Laboratoire de Math\'{e}matiques,
 Universit\'{e} de Franche-Comt\'{e},
 25030 Besan\c con, cedex - France\\
 ericard@univ-fcomte.fr -- AMS classification 46L53}

The classical theory of semigroups has many applications and connections 
with ergodic theory martingales and probability (see \cite{St}). 
The recent developments of the
 non commutative integration in von Neumann provide analogues of these 
notions (\cite{AD, JLMX,JX}). 
For instance, classical Markov semigroups on a probability space
are generalized to  semigroups of unital completely positive maps preserving
a given faithful state. It is natural to try to adapt techniques from the 
commutative theory to the non commutative one. Dealing with $C^*$-algebras,
Sauvageot \cite{Sa} has given a construction of Markov 
$C^*$-dilation for a semigroup
in the spirit of Daniell-Kolmogorov. One of the main tools in the classical 
setting is Rota's dilation theorem \cite{Ro,St},
 it states that any Markovian map 
(unital, positive, self-adjoint on $L_2$ 
and contractive on all $L_p$'s) has a nice
dilation in terms of a reversed martingale, namely $Q^{2n}=\hat 
E\circ E_n$, where $E_n$ are conditional expectations from a decreasing 
filtration and $\hat E$ another conditional expectation. This is closely
related to the construction of Markov chains. Anantharaman-Delaroche
states a counterpart of it for von Neumann algebras in \cite{AD},
 unfortunately some extra
technical condition is needed, she called it ``factorization'' or Markov 
dilation. It is unknown if this condition holds for any Markov operator.
 The aim of 
this note is to discuss this factorization for some concrete and basic 
examples.

\mk

 Let's start with  Stinespring's dilation theorem for $C^*$-algebras. It 
asserts that given a 
unital completely positive map $u:A\to \B(H)$, one can find a Hilbert space 
$K$ containing $H$, a representation $\pi$ from $A$ to $\B(K)$ 
so  that $u$ is just the 
composition of $\pi$ and the natural conditional expectation onto
$\B(H)$. More precisely, $K$ is the Hilbert space 
$A\tens_u H$ obtained
by completion of $A\tens H$ for the scalar product :
$$\forall\; a,a'\in A,\; h,h'\in H \qquad 
\la a\tens h , a'\tens h'\ra_K=\la h, u(a^*a').h'\ra_H$$
The inclusion map $\pi:A\to \B(K)$ is given by $\pi(a)(a'\tens h)=aa'\tens h$.
The embedding  from $H$ into $K$ is  $h\mapsto 1\tens h$. 
The formula for the projection $P$ from $K$ onto $H$
is $P(a\tens h)=1\tens u(a).h$.
And one has 
$$u(a)=P \pi(a)_{|H}.$$
This result is used as the very basic step in Sauvageot's construction.
Its main drawback is that when dealing
with von Neumann algebras (with faithful state or trace) for
applications in $L_p$ spaces (see \cite{AD,JLMX, JX}),
 one would like to have a dilation that 
stays in this category. Therefore Claire Anantharaman-Delaroche introduces 
the notion 
of factorizable maps in \cite{AD}, that we describe precisely now. We will use 
classical notations about von Neumann algebras like in \cite{T,AD}.

\sk

 Let $M$ and $N$  be  von Neumann algebras with normal 
faithful states $\phi$ and $\psi$. 
\begin{defi}
A $(\phi,\psi)$-Markov operator 
$u:(M,\phi)\to (N,\psi)$ is a normal 
unital completely positive
map so that $\psi\circ u=\phi$ and that intertwines the modular groups
 of $\phi$ and $\psi$ ($u\circ \sigma_t^\phi=
\sigma_t^\psi \circ u$).

 One says that $u$ admits a factorization 
if there exist another von Neumann algebra $\tilde M$ with a faithful state
$\tilde \phi$ and normal representations
$\pi : M \to \tilde M$, $\rho : N\to \tilde M$ that are $(\phi,\tilde \phi)$- 
and $(\psi,\tilde \phi)$-Markov maps, with
$$\psi (u(m)n)=\tilde\phi (\pi(m)\rho(n)).$$
\end{defi}
 We say that $(\tilde M,\tilde \phi)$ is a Markov dilation
for $u$.  The conditions on modular groups imply that there is actually a 
$(\tilde \phi,\psi)$-Markov conditional expectation $\E: \tilde M \to N$, and 
$u=\E\circ \pi$. 

 A natural question is to know which maps are factorizable. The aim of this 
note is to give a positive answer for real multipliers.

Viewing $N$ acting as a sub-algebra of $\B(L_2(N,\phi))$, one can 
notice that if $u$ is factorizable then $M\tens_u L_2(N,\phi)$ consists 
exactly in the norm closure of $\pi(M).\rho(N)$ in $L_2(\tilde M,\tilde \phi)$.
  
 In the commutative setting, Stinespring's dilation is actually a Markov 
dilation. This follows from the fact that $N$ acts on $M\tens_u L_2(N,\phi)$
by right multiplications. The commutative von Neumann generated by $M$ and $N$
in $\B(M\tens_u L_2(N,\phi))$ is $\tilde M$ and the state is 
$\tilde \phi (x)=\la 1\tens 1, x.(1\tens 1)\ra$. This is the classical Markov 
construction.

 In the sequel, we are only interested in the case $M=N$ and $\phi=\psi$, and 
we keep the above notations.

We start with some remarks on the set of factorizable maps.

Any state preserving homomorphism $\alpha : M\to M$ is factorizable, 
a dilation
is obtained with $\tilde M=M$, $\tilde \phi=\phi$, $\pi=\alpha$ and $\rho=\Id$.

Any markovian conditional expectation $\E: M\to M$ onto a sub-algebra
$N$ is factorizable. A dilation is given by the free product with
amalgamation over $N$ : $(\tilde M, \tilde \phi)=(M,\phi)*_N(M,\phi)$,
$\pi$ is the homomorphism onto the first copy of $M$ in $\tilde M$ and $\rho$
onto the second one. We refer to \cite{VDN,RX} for definitions.
In particular taking $N=\C. 1$, $\phi$ is factorizable.

The dilation is non unique in general. For instance, $\phi$
 can also be dilated in 
$(M,\phi)\overline \tens_{\min}(M,\phi)$ with the obvious inclusions. 

 The fixed point algebra $N$ by a $\phi$-Markov map plays a particular role
in the 
dilation as we must have $\pi(n)=\rho(n)$ for any $n\in N$. 

If $T$ is $\phi$-Markovian, then it is known that the adjoint of $T$ on 
$L^2(M,\phi)$ comes also from a $\phi$-markovian map denoted by $T^\star$(see 
\cite{AD}),
 that is 
$$\phi(xT(y))=\phi(T^\star(x)y).$$ 
 We say that $T$ is self-adjoint if $T^\star=T$.

\begin{prop}
The set of factorizable $\phi$ Markov operators on $M$ is convex, 
stable by composition, by the involution $^\star$ and 
closed for the point weak-* topology.
\end{prop} 

\begin{pf}
About the involution, it suffices to exchange the role of $\pi$ and $\rho$ as
for analytic elements 
\begin{eqnarray*}
\phi(T^\star(y)x)&=&\phi(T(\sigma_{i}^\phi(x))y)=\tilde\phi(
\pi(\sigma_{i}^\phi(x))\rho(y))\\ &=&
\tilde\phi(\sigma_{i}^{\tilde\phi}(\pi(x))\rho(y))=\tilde\phi(\rho(y)\pi(x))
\end{eqnarray*}

Let $u_i$ be $\phi$-markovian on $M$ with dilation $(\tilde M_i,\tilde \phi_i)$
and morphisms $\pi_i$ and $\rho_i$.

A dilation for $\lambda u_1+ (1-\lambda) u_2$ is given by 
$$(\tilde M_1 \oplus \tilde M_1, \lambda \tilde \phi_1 \oplus 
(1-\lambda) \tilde \phi_2)$$
with morphisms $M\to \tilde M_1\oplus \tilde M_2$, $\pi_1\oplus \pi_2$ and  
$\rho_1\oplus \rho_2$.

A dilation for $u_2\circ u_1$ comes from the free product construction. 
Consider $(\tilde M,\tilde \phi)$ given by 
$(\tilde M_1,\tilde \phi_1)*_M(\tilde M_2,\tilde \phi_2)$ where the 
amalgamation is taken over the copy of $M$ coming from $\rho_1$ in $M_1$ and
$\pi_2$ in $M_2$ 
(note that there are indeed conditional expectation onto them). Let $\E$ be 
the conditional expectation onto the amalgamated copy of $M$. 
From the definition of a dilation,  we get that 
$\E(\pi_1(x))=\rho_1(u_1(x))$ 
and $\E(\rho_2(y))=\pi_2(u_2^\star(y))$ for $x,y\in M$.
A classical computation in free products gives 
\begin{eqnarray*}
\tilde\phi(\pi_1(x)\rho_2(y))&=&\tilde\phi(\E(\pi_1(x))\E(\rho_2(y)))=
\tilde\phi(\rho_1(u_1(x))\pi_2(u_2^\star(y)))\\
&=&\tilde\phi(\rho_1(u_1(x)u_2^\star(y)))=
\phi((u_2\circ u_1)(x).y)
\end{eqnarray*}

 The statement about the closure property is obtained by taking an ultraproduct
(see \cite{R}) cutting with some projections to make representations
normal and the state faithful, technical details can be found in 
\cite{J}.
\end{pf}

 Among other permanence properties, it was observed in \cite{JLMX},
 that the free 
product of factorizable maps is still factorizable, the dilation is simply the
free product of the dilations.

As a corollary of Proposition 2,
 it was pointed to us by Claire Anantharaman-Delaroche that 
any tracial Markov map on $\M_2$ is factorizable as the extreme points of
such maps are exactly the automorphisms.

We now come to Schur multipliers on $\B(\ell_2^I)$  with canonical 
orthonormal basis $(e_i)_{i\in I}$. We will assume that the state $\phi$
has a diagonal density $D=\sum \lambda_i e_i\tens e_i$ for the canonical 
basis with respect to the trace. We have
$\lambda_i>0$ and $\sum \lambda_i=1$.
The modular group of $\phi$ is $\sigma_t^\phi(x)=D^{-it}xD^{it}$.

We represent elements in $\B(\ell_2^I)$ as 
matrices $\M_I$ of size $I$. Given any matrix $T=(t_{i,j})\in \M_I$,
 we say that
$T$ is a Schur multiplier if the following map is well defined
$$M_T=\left\{\begin{array}{ccc}
\B(\ell_2^I)&\to&\B(\ell_2^I)\\
(x_{i,j})&\mapsto& (t_{i,j}x_{i,j})
\end{array}\right.$$
 A characterization of bounded multipliers can be found in \cite{P}. A
 multiplier is (completely) positive
if and only if its symbol $T$ is positive in the sense that for any finite 
set $F\subset I$, $(t_{i,j})_{i,j\in F}$ is positive in $\M_F$. $M_T$ is unital
if $t_{i,i}=1$ for all $i\in I$. $M_T$ is normal unital and completely positive
iff there exist norm 1 vectors $x_i\in \ell_2^I$ so 
that $t_{i,j}=\la x_i,x_j\ra$.   

From these observations, it is clear that 
any unital completely positive Schur multiplier is $\phi$-markovian.

 In the opposite way, if a map $u$ commutes with the modular group of $\phi$
 and $(\log(\lambda_i))_{i\in I}$ is 
independent over $\mathbb Q$ in $\R$, then $u$ has to be a Schur multiplier.

The adjoint of $M_T$ is $M_{T^t}$ where $T^t=(t_{j,i})_{i,j}$. So any 
self-adjoint $\phi$-Markov Schur multiplier has to have real coefficients.

\begin{thm}
Any positive self-adjoint $\phi$-markovian Schur multiplier is factorizable.
\end{thm}

 To construct the dilation, we will need the fermion algebras. Let
$K$ be a real Hilbert space with complexification $K_\C$.
 We briefly recall their construction
and the more general $q$-deformed algebras in the spirit of
$\cite{BKS}$.  The $q$-Fock ($-1\leq q<1$) space over $K$ is
$$\mathcal F_{q}(K)=\C \Omega \oplus \bigoplus_{k} K_\C^{\tens_q k}$$
Where the scalar product on $K_\C^{\tens_q n}$ is given by
$$\la k_1\tens...\tens k_n, h_1\tens ...\tens h_n\ra_q=
\sum_{\sigma\in S_n} q^{|\sigma|} \la k_i,h_{\sigma_i}\ra_{K_\C}$$ 
where $S_n$ is the symmetric group  and $|\sigma|$ the number of inversion 
of the permutation $\sigma$.

When $q=-1$, this is just the antisymmetric tensor product $K_\C^{\wedge n}$.

The creation operator for $e\in K$ is given by
$$l(e).( h_1\tens...\tens h_n) = e \tens h_1\tens...\tens h_n$$
They satisfy the $q$-relation
$$l(f)^*l(e)-ql(e)l(f)^*=\la f,e\ra_K Id$$
The $q$-von Neumann algebra is
$$\Gamma_{q}(K)=\{ \omega(e)=l(e)+l(e)^* \; ; \; e\in K\}''$$
It is type $II_1$ with the trace $\tau (x)=\la \Omega, x.\Omega\ra_{
\mathcal F_q(K)}$.

We are mainly concerned with the fermion algebra when $q=-1$. If $e\in K$
has norm 1, then $\omega(e)$ is a symmetry, i.e. self adjoint with 
$\omega(e)^2=1$. Moreover we have 
$$\forall e,f\in K \qquad \tau(\omega(e)\omega(f))=\la e,f\ra_{K}$$

\begin{pf}
We will use the notation $e_{i,j}$ for the canonical basis of $\B(\ell_2^I)$.

Let $M_T$ be a self-adjoint markovian Schur multiplier. As 
$t_{i,j}=\la x_i,x_j\ra$ is real, $T$ defines a new 
scalar product on the real linear span of $e_i$'s by the formula :
$$\la \sum a_i e_i,\sum b_i e_i\ra_T= 
\la \sum a_i x_i,\sum b_i 
x_i\ra_{\ell_2^I}=\sum_{i,j}  a_i b_j t_{i,j}$$
We call $\ell_{2,T}$ the real 
Hilbert space obtained after quotient and completion.
We still denote by $e_i$ the representative of $e_i$ in $\ell_{2,T}$, we have
$$\la e_i,e_j\ra_T=t_{i,j}=t_{j,i}$$
 
 Let $\tilde M=\B(\ell_2^I)\overline \tens_{\min} \Gamma_{-1}(\ell_{2,T})$ with
normal faithful state $\tilde \phi=\phi\tens \tau$.  Let 
$$d=\sum_i e_{i,i} \tens \omega(e_i)\in \tilde M$$
It is a unitary (symmetry) in the centralizer of $\tilde \phi$, as 
$e_{i,i}$ are in the one of $\phi$.

For $x\in \tilde M$, define 
$$\mathcal U(x)=dxd$$
This is a $\tilde \phi$ Markov map, and a representation as $d^2=1$.

Let $\pi: M\to \tilde M$ be the obvious inclusion $\Id\tens 1$. Define 
$\rho : M\to \tilde M$ as $\rho=\mathcal U\circ \pi$. It is clear that 
$\pi$ and $\rho$ are $(\phi,\tilde \phi)$-markovian and for $x=(x_{i,j})$
and $y=(y_{i,j})$ finite matrices : 
\begin{eqnarray*}
\tilde \phi (\pi(x)\rho(y))&=&\phi\tens \tau 
\Big((x\tens 1) d (y\tens 1)d \Big)\\
&=& \phi\tens \tau \Big( (x_{i,j}.1). (y_{i,j}\omega(e_i)\omega(e_j))\Big)\\
&=&\sum_{i,j} \lambda_i x_{i,j}y_{j,i} \tau (\omega(e_i)\omega(e_j))\\
&=&\sum_{i,j} \lambda_i x_{i,j}y_{j,i} t_{i,j}\\
&=& \phi (T(x)y)
\end{eqnarray*}
\end{pf}

When $T=\Id$, $M_T$ is a conditional expectation and this dilation is very 
different from the one obtained by free product.

 Combining the previous example with the permanence properties gives a 
wide class of factorizable maps. If $I$ is finite and $\phi$ is the trace, we
can take compositions of multipliers in different basis (and with conditional 
expectations, representations) and convex 
combinations of them. We do not know whether we can achieve all tracial 
markovian maps for $\M_n$. It follows from Grothendieck's theorem that any 
completely bounded Schur multiplier is a multiple (less than the Grothendieck
constant) of a convex combination of rank one multipliers (see \cite{P}). 
In terms of Markov maps, any Markov Schur 
multiplier can be obtained as a linear combination of representations, 
unfortunately it can not be a convex one (or the Grothendieck constant would 
have to be exactly 1).

\sk

 If one looks carefully, the von Neumann  algebra generated by 
$\pi(\B(\ell_2^I))$ and $\rho(\B(\ell_2^I))$ is exactly $\B(\ell_2^I)\overline
\tens \Gamma_{-1}^{e}(\ell_{2,T})$ where $\Gamma_{-1}^{e}(\ell_{2,T})$ 
is the sub-algebra of $\Gamma_{-1}(\ell_{2,T})$ generated by even elements of
the form $\omega(e_i)\omega(e_j)$.

 Now that we have a dilation for Schur multipliers, thanks to the construction
in \cite{AD} chapter 6, we have a non commutative Markov chain. In our concrete
setting, one can avoid this abstract construction (and free products). 
We follow the notations of \cite{AD} and classical ones for infinite tensor 
products (we drop the completion symbol),
let $$M=\B(\ell_2^I)  \tens
\Gamma_{-1}^{e}(\ell_{2,T})^{  \tens \infty}\subset 
\Gamma_{-1}(\ell_{2,T})^{  \tens \infty}$$ 
equipped with the tensor product
state. Actually because of the commutation relations
$$\Gamma_{-1}^{e}(\ell_{2,T})^{  \tens \infty}\subset 
\Gamma_{-1}(\ell_{2,T}\tens \ell_2)$$
 Let $J_0:\B(\ell_2^T)\to M$ be the natural inclusion given by 
$$J_0(x)=x\tens 1\tens 1
 \tens ...$$
 The letter $S$ stands for the shift on 
$\Gamma_{-1}^{e}(\ell_{2,T})^{  \tens \infty}$:
$$S(x_1\tens...\tens x_n\tens 1\tens  ...)=1\tens x_1\tens...\tens x_n\tens 1
\tens  ...$$
We will also need the symmetry 
$$d_1=\sum_i e_{i,i} \tens \omega(e_i)\tens 1 \tens ... \in 
\B(\ell_2^I)\overline
\tens \Gamma_{-1}(\ell_{2,T})^\infty$$ We have an injective 
morphism $\beta:M \to M$ given by 
$$\beta(x)=d_1S(x)d_1$$
If  $J_q=\beta^q\circ J_0$, the $q^{th}$ copy of $\B(\ell_2^I)$ is 
$$J_q(e_{i,j})= e_{i,j} \tens\underbrace{ \omega(e_i)\omega(e_j) \tens
...\tens \omega(e_i)\omega(e_j)}_{q\textrm{ times}} \tens 1 \tens ...$$
The algebra $\mathcal B_{n]}$ generated by the first $n^{th}$ copies of 
$\B(\ell_2^I)$ is exactly 
  $$\mathcal B_{n]}=\B(\ell_2^I)  \tens
\Gamma_{-1}^{e}(\ell_{2,T})^{  \tens n} \tens \C^{\tens\infty}$$
And $\mathcal B_{[n}$, generated by $J_q(\B(\ell_2^I))$ with $q\geq n$, is 
$$\mathcal B_{[n}=J_n(\B(\ell_2^I)). ( \C^{\tens n} \tens
 \Gamma_{-1}^{e}(\ell_{2,T})^{  \tens \infty})$$
All maps preserve the involved modular groups and if $\E_{n]}$ and $\E_{[n}$ 
are the conditional expectations onto $\mathcal B_{n]}$ and 
$\mathcal B_{[n}$, one can 
check the Markov properties 
\begin{eqnarray*}
\E_{n]}\circ J_q&=&J_n \circ T^{q-n} \qquad q\geq n\\
\E_{n+q]}\circ \beta^q&=&\beta^q\circ \E_{n]}\\
\E_{[n}\circ J_0&=&J_n\circ T^n
\end{eqnarray*}
In particular Rota's dilation is 
$$J_0\circ T^{2n}=\E_{0]} \circ\E_{[n} \circ J_0$$

\mk

The above construction can also be carried out for Fourier multipliers on 
discrete groups. Let $G$ be a discrete countable group and $L(G)\subset 
\B(\ell_2G)$ its left von Neumann algebra. It is the bicommutant of
the left translations by $g\in G$, denoted as usual by $\lambda(g)$. 
It is a type $II_1$ algebra with trace given by 
$$\tau(\lambda(g))=\la \delta_e,\lambda(g)\delta_e\ra=\delta_{g,e}$$  
 A function  $t:G\to \C$ defines a Fourier multiplier $M_t$ 
if the following map is well defined on $L(G)$ 
$$M_t(\lambda(g))= t_g\lambda(g)$$
Actually a Fourier multiplier is completely bounded iff  the Schur
multiplier $(t_{h^{-1}g})_{h,g}$ is bounded on $\B(\ell_2G)$ (see \cite{P}). 
It is unital 
completely positive iff $t_e=1$ ($e$ is the unit of $G$) 
and $t$ is positive definite, in this
case it is $\tau$ markovian.
$M_t$ is self-adjoint if $t_g=t_{g^{-1}}\in \R$.

\begin{cor}
Any self-adjoint unital completely positive Fourier multiplier on a discrete
group is factorizable.
\end{cor}

\begin{pf}
The dilation can be obtained directly from the previous one but can be 
reinterpreted in terms of a cross product as follows.

Let $\ell_{2,T}$ be the Hilbert space obtained by quotient and completion 
of the real span of $\lambda(g)$ for the scalar product :
$$\la \lambda(g),\lambda(h)\ra_T=t_{g^{-1}h}$$

We let $h$ be the class of $\lambda(h)$ in $\ell_{2,T}$.
 It is clear that $G$ acts unitarly on $\ell_{2,T}$ by left multiplications.
So $G$ also acts by automorphisms on $\Gamma_{-1}(\ell_{2,T})$ by, for 
$g,h\in G$ 
$$\alpha(g).\omega(h)=\omega(gh)$$
Let $\tilde M$ be the crossed product 
$\Gamma_{-1}(\ell_{2,T})\rtimes_\alpha G$ (see \cite{T} for definitions).
 This is again a type $II_1$ 
algebra, $\tilde \phi$ is the canonical trace and for $g\in G$, $x\in 
\ell_{2,T}$
$$\lambda(g) \omega(x)\lambda(g^{-1})=\alpha(g).\omega(x)=\omega(g.x)$$  

Then $\pi$ is just the natural copy of $L(G)$ in $\tilde M$ and $\rho$ is 
given by $\rho(x)=\omega(e)\pi(x)\omega(e)$ for $x\in L(G)$.
 Note that $\omega(e)$ is a symmetry. 
We have
\begin{eqnarray*}
\tilde \phi (\pi(\lambda(g))\rho(\lambda(h)))
&=&\tilde\phi ( \lambda(g)\omega(e)\lambda(h)\omega(e))\\
&=&\tilde\phi ( \lambda(g)\lambda(h)\omega(h^{-1})\omega(e))\\
&=&\delta_{gh,e}t_h=\delta_{gh,e}t_g\\
&=&\tau (M_t(\lambda(g)) \lambda(h))
\end{eqnarray*} 
\end{pf}

As before, one has nice formulas for the associated Markov chain.
Briefly speaking, 
$M=\Gamma_{-1} (\ell_{2,T}\tens \ell_2) \rtimes_\alpha G$, where
$\alpha$ is the diagonal action of $G$ : $\alpha(g).
\omega(h\tens v)=\omega(gh \tens v)$ for $g,h\in G$ and $v\in\ell_2$.
 The inclusion $J_0$ of $L(G)$ is the 
natural one. The morphism $\beta$ is given by the shift 
on the canonical basis $(e_i)$ of $\ell_2$ and a conjugation 
$$\beta (\lambda(g)\omega(h\tens e_i))=\omega(e\tens e_0)\lambda(g)
\omega(h\tens e_{i+1}) \omega(e\tens e_0)$$
 Then one defines 
$\mathcal B_{n]}$ and 
$\mathcal B_{[n}$ as above to get Rota's construction.

\mk

There is another family of maps close to multipliers which can be seen to have
a dilation quite easily. It consists of the maps arising from the second 
quantization on $q$-deformed algebras (see \cite{BKS}).
 For any real Hilbert space $K$, let 
$\omega_K: \Gamma_q(K)\to \mathcal F_q(K)$ be $\omega_K(x)=x.\Omega_K$. This
 is an 
injective mapping with dense range, we denote by
 $\omega_K^{-1}$ its inverse (not defined everywhere).
 
If $K$ and $L$ are real Hilbert spaces, any contraction $T:K\to L$ gives 
rise to a tracial Markov map $\Gamma(T):\Gamma_q(K)\to \Gamma_q(L)$ 
satisfying :
$$T(\omega_K^{-1}(k_1\tens...\tens k_n))=
\omega_L^{-1}(T(k_1)\tens...\tens T(k_n))$$

If $T$ is isometric (unitary) then, 
$\Gamma(T)$ is an injective representation (automorphism). 
So if $K\subset L$, we can see $\Gamma_q(K)$ as a sub-algebra of $\Gamma_q(L)$.
In this situation, if $P$ is the orthogonal projection of $L$ onto $K$,
$\Gamma(P)$ is then the trace preserving 
conditional expectation from $\Gamma_q(L)$ onto $\Gamma_q(K)$.

Any contraction $T: K\to K$ 
can be dilated to a unitary say $U$ (to a symmetry if $T$ is self-adjoint)
 on a bigger Hilbert space $L$. From $T=PU_{|K}$, one sees that a dilation 
of $\Gamma(T)$ is given by $\Gamma_q(L)$ with its trace, with 
$\pi$ the natural injection and $\rho=\Gamma(U^*)\circ \pi$.  

 To get Rota's dilation is also easy in this case, assuming that $T$ 
is self-adjoint.
 Indeed, let $U$ be a 
strong dilation of $T$ on $L$, that is $PU^k_{|K}=T^k$ (see \cite{P} 
Theorem 1.1). Let $K_n={\rm span}\{ U^l(K)\,;\, l\geq n\}$ and $P_n$ the 
corresponding projection, then one has 
$$P P_nP = T^{2n}$$
Indeed if $k\in K$ and $x\in K$ for $l\geq 0$
$$\la k, U^{n+l}(x)\ra=\la k,P U^{n+l}P(x)\ra = \la T^{n+l}(k),x\ra
=\la T^n(k),U^{l}(x)\ra=\la U^n(T^n(k)),U^{n+l}(x)\ra$$
so $P_nP(k)=U^{n}(T^n(x))$.

Going to second quantization, with $\E_n=\Gamma(P_n)$
the conditional expectation from $\Gamma_{q}(L)$ 
onto $\Gamma_{q}(K_n)$, and $\hat \E=\Gamma(P)$ being the
  conditional expectation  onto $\Gamma_{q}(K)$ (adjoint of the inclusion 
$J$), one has
$$\hat \E \circ \E_n \circ J = \Gamma(PP_nP)\circ J=J\circ \Gamma(T)^{2n} $$

 The same dilation works also for the $q$-deformed versions of Araki-Woods 
factors of Hiai (\cite{H}).

\bibliographystyle{plain}

\end{document}